\documentclass{article}
\def \Z {{\mathbf {Z}}}
\def \T {{\mathbf {T}}}
\def \X {{\mathbf {X}}}
\def \N {\mathbf N}
\def \S {{\mathbf {S}}}

\def\uu{\bigsqcup}

\usepackage[T1]{fontenc} 
\usepackage[cp1251]{inputenc} 
\usepackage[russian]{babel}

\title{ Некоммутирующие автоморфизмы с расходящимися эргодическими средними}

\author{ В.В. Рыжиков}
\date{}

\begin{document}
\Large

\maketitle

\begin{abstract}  Остин и независимо Хуанг, Шао и  Е нашли примеры пар некоммутирующих автоморфизмов, обладающих расходящимися двукратными эргодическими  средними.
Мы предлагаем новые примеры среди гауссовских и пуассоновских надстроек 
над  специально подобранными преобразованиями с инвариантной сигма-конечной мерой.

\vspace{2mm}
\it Ключевые слова:  \rm расходимость эргодических средних,  преобразования ранга один, гауссовские и пуассоновские надстройки, нулевая энтропия.  
\end{abstract}

\section{Введение}
 В статье \cite{F} был задан вопрос о сходимости в пространстве $ L_2(X, \mu)$   средних  для последовательности вида  
$\{S^nf\, T^ng\}$, где $f,g\in L_\infty(\mu)$, а  $S,T$ -- детерминированные автоморфизмы вероятностного пространства $(X,\mu)$.  Отметим, что авторы \cite{F} доказали, что для таких $S,T$   при $p\geq 2$ имеет место сходимость средних вида  
$$ \frac 1 N \sum_{n=1}^N S^nf\,T^{n^p}g.$$
В  работах  \cite{A} и  \cite{Ye} разными способами построены  подходящие  пары  $T,S$ и  $f,g\in L_\infty(X, \mu)$ с расходящейся  последовательностью
$$ \frac 1 N \sum_{n=1}^N \int T^nf\,S^ng \, d\mu.$$

В нашей заметке предлагается  еще одно решение, основанное на применении  преобразований  с сигма-конечной инвариантной мерой, будем называть их бесконечными преобразованиями.  Пуассоновские и гауссовские надстройки над  подходящими преобразованиями  дают упомянутый эффект расходимости. 
Опишем свойства наших автоморфизмов  $S,T$. Для некоторого множества $C$,  на большом временном промежутке для времени $n$ из этого промежутка выполняется   $S^nC=T^nC$, а потом для $n$ из
другого более долгого промежутка  множества $S^nC$  и $T^nC$ являются независимыми. Эту ситуацию легко реализовать для бернуллиевского сдвига $S$ и сопряженного с ним $T$.  А  для реализации с нулевой энтропией роль таких $S$, $T$ успешно   играют  пуассоновские  надстройки. В нашем случае $S=S'\times S'$, множество $C$  и все  $S^n C$ принадлежит координатной алгебре первого сомножителя $S'$.  Степени  $T^n$ таковы, что    длительное время, начиная с какого-то момента, мы наблюдаем $T^nC=S^nC$, потом $T$-динамика начинает постепенно перекачивать множества $T^nC$  в из одной в другую координатную алгебру. Начинается  временной интервал, когда множества    $T^nC$ долго пребывают во второй координатной  алгебре автоморфизма $S$ и по этой причине независимы  от множеств $S^nC$, лежащих в первой координатной алгебре.  Потом происходит перекачка обратно и наступает долгий период, когда  $T^nC=S^nC$. 
Далее возрастающие интервалы  совпадений  чередуются с интервалами  независимости множеств $T^nC$ и  $S^nC$, обеспечивая расходимость 
обсуждаемых  средних.

\vspace{3mm}
\bf Теорема 1. \it Существуют бесконечные автоморфизмы $\S,\T$,  множество $A$ положительной конечной меры
и последовательность $\{N_m\}$ такая, что $|N_{m+1}|/|N_m|\to\infty$,  
причем $T^iA\cap S^iA=\phi$ для $i\in [N_{4k}, N_{4k+1}]$ и
 $T^iA=S^iA$ для $i\in [N_{4k+2}, N_{2k+3}]$.

Для соответствующих гауссовских  и пуассоновских надстроек $S,T$ над $\S,\T$ для некоторой неотрицательной функции  $f\in L_\infty$  будет расходиться последовательность
$$ a_N=\frac 1 N \sum_{n=1}^N \int T^nf\,S^nf\, d\mu.$$
При этом предельные точки последовательности $a_N$ образуют отрезок $[c^2, c]$,\ $c=\int f\, d\mu$. 
Искомые  $S,T$ могут быть как неперемешивающими, так и  перемешивающими автоморфизмами с  нулевой энтропией.\rm

\vspace{3mm}

\section{\ Бесконечные   преобразования  $\bf \T,\S$.} 
Приступим к построению нужных примеров бесконечных преобразований.
Ищем их  в виде  $\Z_2$-расширений некоторой  конструкции  $R$ ранга один.
Определим ее.  

Пусть заданы параметры 
   $h_1=1$, $r_j\geq 2$ и наборы натуральных  чисел
$$ \bar s_j=(s_j(1), s_j(2),\dots,s_j(r_j)), \ s_j(i)\geq 0, \ j\in\N. $$
Фазовое пространство $X$ для преобразования $R$, определенного ниже,  представляет собой объединение башнен $$X_j=\uu_{i=0}^{h_j-1} R^iE_j,$$
где  $R^iE_j$ — непересекающиеся полуинтервалы, называемые этажами. 

На этапе $j$ преобразование  $R$ определено как обычный перенос интервалов (этажей), но на верхнем этаже 
$R^{h_j-1}E_j$ башни $X_j$ преобразование пока не определено.
 Башня $X_j$ разрезается на $r_j$ одинаковых узких подбашен $X_{j,i}$ (они называются колоннами),    и над каждой колонной  $X_{j,i}$  добавляется  $s_j(i)$ новых этажей. Преобразование $R$ по-прежнему определяется как подъем на этаж выше, но  самый последний надстроенный этаж над колонной с номером  $i$ преобразование $R$ отправляет в нижний этаж колонны   с номером  $i+1$.  Таким обазом возникает   новая башня 
$$X_{j+1}=\uu_{i=0}^{h_{j+1}-1} R^iE_{j+1}$$
высоты 
   $$h_{j+1}=r_jh_j+\sum_{i=1}^{r_j}s_j(i),$$
где $E_{j+1}$ -- нижний этаж колонны $X_{j,1}$.
Отметим, что доопределяя преобразование, мы полностью сохраняем предыдущие построениея. 
  Продолжая этот процесс до бесконечности,  
получаем обратимое преобразование $R:X\to X$, сохраняющее меру  Лебега на объединении $X=\uu_j X_j$.

\vspace{3mm}
\bf Автоморфизмы $\bf \T,\S$. \rm
Пусть $R$ — конструкция ранга один с параметрами  $r_j=j$, $s_j(i)=jh_j$.
Положим $$\X = X\times \Z_2, \ \ \
E=\uu_j (R^{h_{2j}}E_{2j}\uu R^{2jh_{2j}}E_{2j}).$$
Определим  $ \T,\S$ следующим образом:

$\S(x,z)=(Rx, z), $

$\T(x,z)=(Rx, z+1), \ \ x\in E,$ 

$\T(x,z)=(Rx, z), \ \ x\notin E.$

\vspace{3mm}
\bf Лемма 2. \it Для  $A=X_1\times {0}\subset\X$
множества    $\T^iA$ и $\S^iA$ не пересекаются при $h_{2j}< i<2jh_{2j}$ и
 совпадают  при  $h_{2j+1}<i< 2jh_{2j+1}$. \rm

\vspace{3mm}
Доказательство. Назовем координату $z$ уровнем. Пусть $x\in X_1$. Заметим, что уровень точек $\T^i(x,z)$ меняется, когда точка $R^ix$ окажется в множестве $R^{h_{2j}}E_{2j}$.  Поэтому $\S^iA$ и $\T^iA$ расположены  на разных уровнях при $h_{2j}< i<2jh_{2j}$. Уровень  снова изменится, когда   $R^ix$ встретит множество $R^{2jh_{2j}}E_{2j}$.  При $h_{2j+1}<i< 2jh_{2j+1}$ выполняется  $\T^i(x,z)=\S^i(x,z)$ и тем самым $\S^iA=\T^iA$.

\section{ Пуассоновские надстройки}
 Рассмотрим конфигурационное пространство $X_\circ$, 
состоящее из всех бесконечных счетных множеств $x_\circ$ таких, что каждый  интервал из пространств $X$ содержит лишь конечное число элементов множества $x_\circ$.

Пространство $X_\circ$ оснащается  мерой Пуассона. Напомним ее определение.
Подмножествам $A\subset X$ конечной $\mu$-меры в конфигурационном пространстве $X_\circ$  сопоставлены цилиндрические  множества
   $C(A,k)$, $k=0,1,2,\dots$,    по формуле
$$C(A,k)=\{x_\circ\in X_\circ \ : \ |x_\circ\cap A|=k\}.$$

Всевозможные конечные пересечения вида $\cap_{i=1}^N C(A_i,k_i)$
     образуют полукольцо.
На этом полукольце определена мера  $\mu_\circ$ следующим образом:
   при условии, что измеримые множества $A_1, A_2,\dots, A_N$ не пересекаются
и имеют конечную меру, положим
$$\mu_\circ(\bigcap_{i=1}^N C(A_i,k_i))=\prod_{i=1}^N \frac {\mu(A_i)^{k_i}}{k_i!} e ^{-\mu(A_i)}.\eqno (\ast)$$
Пояснение: если множества $A$, $B$ не пересекаются, то
    вероятность $\mu_\circ(C(A,k))\cap C(B,m))$  одновременного появления $k$ точек конфигурации $x_\circ$ в $A$
     и $m$ точек конфигурации $x_\circ$ в $B$ равна произведению вероятностей $\mu_\circ(C(A,k))$ и $\mu_\circ( C(B,m))$.
    Другими словами, события $C(A,k))$ и $C(B,m)$ независимы. Так как  множества $A_1, A_2,\dots,A_N$ не пересекаются,     в формуле $(\ast)$ фигурирует произведение соответствующих вероятностей. 
 
Классическое продолжение меры является  пространством Пуассона
$(X_\circ,\mu_\circ)$,     изоморфным стандартному вероятностному пространству Лебега.
Автоморфизм $T$ пространства $(X,\mu)$ естественным образом индуцирует автоморфизм
$T_\circ$ пространства $(X_\circ,\mu_\circ)$,  называемый пуассоновской надстройкой над $T$.

\section{Доказательство теоремы 1}
Пусть $A=X_1\times \{0\}\subset\X$,  рассмотрим цилиндр $C=C(A,m)$,
 фиксируя $m$.  Из леммы 2 и определений пуассоновских надстроек  непосредственно получаем следующее утверждение.

\vspace{2mm}
\bf Лемма 3. \it Обозначим $T=\T_\circ,$ $S=\S_\circ$. 
Множества $S^{i}C$ и $T^iC$ независимы при  $h_{2j}< i<2jh_{2j}$ и
 совпадают  при  $h_{2j+1}<i< 2jh_{2j+1}$. \rm

\vspace{2mm}
Пусть $f=\chi_C$, $\mu_\circ(C)=c,$
тогда из леммы 3 получим  
$$ \frac 1 {2jh_{2j}}\sum_{i=1}^{2jh_{2j}}\int S^if\, T^if\ d\mu_\circ \to \ c^2$$
и
 $$ \frac 1 {2jh_{2j+1}}\sum_{i=1}^{2jh_{2j+1}}\int S^if\, T^if\ d\mu_\circ  \ \to \ c.$$
Автоморфизмы  $S$ и $ T$
 имеют нулевую энтропию. Пуассоновская надстройка $S'$ над конструкцией ранга один имеют нулевую энтропию, поэтому 
энтропия $S=S'\times S'$ также  нулевая. 

Автоморфизмы  $S$ и $T$  сопряжены, 
так как  $\bf VSV= T$, где инволюция  $\bf V$ определяется следующим образом: 
${\bf V}(x,z)=(x,z+1)$   для $x$ из множества 
$\uu_j \uu_{i= h_{2j}+1}^{2jh_{2j}} R^iE_{2j}$, а для других  $x$ выполнено 
${\bf V}(x,z)=(x,z)$.  

\bf Перемешивающие примеры $S,T$. \rm Чтобы  получить перемешивающие преобразования имеется много разных возможностей. Воспользуемся, например,    бесконечными конструкциями  $R$ ранга один  из работы \cite{DR}. Пусть  $r_j=j$ и для большинства этапов выполнено  $s_j(i)=jh_j+i$ (лестничные надстройки).
Другие   этапы  пусть  обеспечивают  сингулярность (и даже простоту) спектра соответствующих гауссовских и пуассоновских надстроек, как это было издожено в см.\cite{DR}.  

Полагаем $$\X = X\times \Z_2$$
и для некоторого бесконечного множества $J$ (отвечающего немодифицированным этапам)   определим
 $$E=\uu_{j\in J} (R^{h_{2j}}E_{2j}\uu R^{2jh_{2j}}E_{2j})$$
и повторим определение преобразований $\S,\T$:
$$\S(x,z)=(Rx, z), $$
$$\T(x,z)=(Rx, z+1), \ \ x\in E, \ \ \ \ \T(x,z)=(Rx, z), \ \ x\notin E.$$
Для соответствующих пуассоновских надстроек также выполняется утверждение леммы  3. 

Для гауссовских настроек  аргументация  расходимости  по сути ничем не отличается от пуассоновского случая. 
Рассматрим   $\S,\T$ как 
 ортогональные  операторы в $L_2$, пусть  $S$ и $T$ обозначают теперь их действия в  вещественном пространстве $L_2$, оснащенном   инвариантной гауссовской мерой. Пользуясь тем, что  ортогональность пространств обеспечивает независимость соответствущих цилиндрических множеств, в качестве $f$ выбираем индикатор цилиндра с основанием на одномерном пространстве вектора $\chi_A$.   Так как   $\S,\T$ обладают сингулярным спектром, 
энтропия надстроек нулевая. Отметим, что гауссовские автоморфизмы   использовались в работе  \cite{A}, для них другим способом подбиралсь подходящим образом ортогональные операторы.

Остается лишь заметить, что при $N_{4j}=h_{2j}$,  $N_{4j+1}=2jh_{2j}$,
$N_{4j+2}=h_{2j+1}$,  $N_{4j+3}=2jh_{2j+1}$ выполняется утверждение теоремы 1.

\section{ Замечания} 
 Автоморфизмы $S$ и $T$ оказались   сопряженными. Это не случайно, поясним почему.  Джойнинг $\nu$ автоморфизмов $S$ и $T$ -- это $(S\times T)$-инвариантная мера на $X\times X$ с проекциями $\mu$ на сомножители в $X\times X$.  Пусть $S$ и $T$ дизъюнктны, то есть  они  обладают  только тривиальным   джойнингом  $\nu= \mu\otimes\mu$.  Известно, что из всякой последовательности $N'\to\infty$ можно выбрать подпоследовательность $N_i$  такую, что для некоторого джойнинга $\nu$ выполнено 
$$ a_{N_i}=\frac 1 {N_i} \sum_{n=1}^{N_i} \int T^nf\,S^ng\, d\mu\ \to\ 
\int f\otimes g d\nu.$$
По условию $S$ и $T$ дизъюнктны,  значит, $\nu=\mu\otimes\mu$, поэтому 
$$a_{N_i}\to \int Tf\,d\mu\int g\, d\mu.$$  
В силу произвольного выбора  $N'\to\infty$ получаем, что   исходная последовательность $a_N$ сходится. Таким образом, для контрпримера необходима недизъюнктность  $S$ и $T$, случай изоморфизма $S$ и $T$ здесь оптимален. 

Отметим, что  для  типичной (в смысле категорий Бэра) пары автоморфизмов $S,T$ верно, что $S$ и $T$ дизъюнктны, поэтому  типична сходимость обсуждаемых средних.

Автоморфизмы $S$ и $T$ дизъюнктны, например, если $S$ обладает нулевой, 
а $T$   вполне положительной $P$-энтропией Кириллова-Кушниренко, что показано в  работе  \cite{RT}.  Подходящие   примеры   в \cite{RT} предъявлены тоже в виде   пуассоновских надстроек над преобразованиями ранга один.

\end{document}